\documentclass[12pt]{article}

\usepackage{latexsym,amssymb, amsthm, epsfig}

\usepackage{a4wide} 
\newcommand{{\cao}}{{\c{c}\~{a}o}}  
\newcommand{{\coes}}{{\c{c}\~{o}es}}  
\def\E{{\mathbb E}}  
\def\P{{\mathbb P}}  
\def\R{{\mathbb R}}  
\def\Z{{\mathbb Z}}  
  
\def\N{{\mathbb N}}  
\def\etatnx{{\frac{\eta_t^N(x)}{N}}}  
\def\etanx{{\frac{\eta^N(x)}{N}}}  
\def\etany{{\frac{\eta^N(y)}{N}}}  
\def\etatny{{\frac{\eta_t^N(y)}{N}}}  
  
\def\prova{\noindent{\bf Proof: }}

\def\ron{{\rho^{N}}}  

\def\cL{{\mathcal{L}}}

\def\qsd{{\sc{qsd}}}  
\def\fv{{\sc{fv}}}

\def\({{\Bigl(}}  
\def\){{\Bigr)}}  
\def\tip{{\ensuremath{{\rm type}}}}  
\def\mua{{\ensuremath{\mu_\alpha}}}

\def\reff#1{(\ref{#1})}
\def\one{{\mathbf 1}}  
\def\ind{{\mathbf 1}}

\def\rI {{\rm I}}
\def\tim{{\rm time}}

\def\lab{{\rm label}}
\def\square{\ifmmode\sqr\else{$\sqr$}\fi}
\def\sqr{\vcenter{
         \hrule height.1mm
         \hbox{\vrule width.1mm height2.2mm\kern2.18mm\vrule width.1mm}
         \hrule height.1mm}}                  

\theoremstyle{plain}   
\newtheorem{lema}{Lemma}[section]   
\newtheorem{teo}{Theorem}[section]   
\newtheorem{propo}{Proposition}[section]

\theoremstyle{definition}

\theoremstyle{remark}

\begin{document}


\centerline {\Large Quasi stationary distributions and Fleming-Viot processes}
\centerline{\Large in countable spaces}

\vspace{0.5cm}

\centerline{\bf Pablo A.  Ferrari, Nevena Mari\'{c} } 

\centerline{\emph{Universidade de S\~ao Paulo}}
 
\vspace{0.5cm}

{\large\bf Abstract}  
We consider an irreducible pure jump Markov process with rates $Q=(q(x,y))$ on
$\Lambda\cup\{0\}$ with $\Lambda$ countable and $0$ an absorbing state. A {\em
  quasi-stationary distribution \rm} ({\qsd}) is a probability measure $\nu$ on
$\Lambda$ that satisfies: starting with $\nu$, the conditional distribution at
time $t$, given that at time $t$ the process has not been absorbed, is still
$\nu$. That is, $\nu(x) = \nu P_t(x)/(\sum_{y\in\Lambda}\nu P_t(y))$, with $P_t$
the transition probabilities for the process with rates $Q$.
   
A {\em Fleming-Viot} ({\fv}) process is a system of $N$ particles moving in
$\Lambda$. Each particle moves independently with rates $Q$ until it hits the
absorbing state $0$; but then instantaneously chooses one of the $N-1$ particles
remaining in $\Lambda$ and jumps to its position. Between absorptions each
particle moves with rates $Q$ independently.
  
Under the condition $\alpha:=\sum_x\inf Q(\cdot,x) > \sup Q(\cdot,0):=C$ we
prove existence of {\qsd} for $Q$; uniqueness has been proven by Jacka and
Roberts. When $\alpha>0$ the {\fv} process is ergodic for each $N$. Under
$\alpha>C$ the mean normalized densities of the {\fv} unique stationary measure
converge to the {\qsd} of $Q$, as $N \to \infty$; in this limit the variances
vanish.

\paragraph{Keywords} Quasi stationary distributions. Fleming-Viot process. 


\section{Introduction}

Let $\Lambda$ be a countable set and $Z_t$ be a pure jump regular Markov process
on $\Lambda\cup\{0\}$ with transition rates matrix $Q=(q(x,y))$, transition
probabilities $P_t(x,y)$ and with absorbing state $0$; that is $q(0,x)=0$ for
all $x\in\Lambda$. Assume that the exit rates are uniformly bounded above: $\bar
q:=\sup_{x} \sum_{y\in\{0\}\cup\Lambda\setminus\{x\}} q(x,y)<\infty$, that
$P_t(x,y)>0$ for all $x,y\in\Lambda$ and $t>0$ and that the absorption time is
almost surely finite for any initial state.  The process $Z_t$ is ergodic with a
unique invariant measure $\delta_0$, the measure concentrating mass in the state
$0$.  Let $\mu$ be a probability on $\Lambda$. The law of the process at time
$t$ starting with $\mu$ conditioned to non absorption until time $t$ is given by
\begin{equation}
  \label{u7}
  \varphi^\mu_t (x) = \frac{\sum_{y\in\Lambda} \mu(y)
  P_t(y,x)}{1-\sum_{y\in\Lambda} \mu(y) P_t(y,0)}\,,~~x \in \Lambda.
\end{equation}

A \emph{quasi stationary distribution} ({\qsd}) is a probability measure $\nu$
on $\Lambda$ satisfying $\varphi^\nu_t = \nu$.  Since $P_t$ is honest and
satisfies the forward Kolmogorov equations we can use an equivalent definition
of {\qsd}, according Nair and Pollett \cite{nair}. Namely, a {\qsd} ( and only a
{\qsd}) is a left eigenvector $\nu$ for the restriction of the matrix $Q$ to
$\Lambda$ with eigenvalue $-\sum_{y\in\Lambda} \nu(y) q(y,0)$: $\nu$ must
satisfy the system
\begin{eqnarray}\label{eq:quasi2}   
  \sum_{y \in \Lambda}\nu(y)\,[q(y,x)+q(y,0)\nu(x)]=0 , ~~\forall x \in \Lambda.   
\end{eqnarray} 
(recall $q(x,x) = -\sum_{y\in\Lambda\cup\{0\}\setminus\{x\}} q(x,y)$.)

The \emph{Yaglom limit} for the measure $\mu$ is defined by
\begin{equation}\label{yl}  
  \lim_{t\to \infty}\varphi^\mu_t (y)\,,\quad y\in\Lambda
\end{equation}  
if the limit exists and it is a probability on $\Lambda$.

When $\Lambda$ is finite, Darroch and Seneta (1967) prove that there exists a
unique {\qsd} $\nu$ for $Q$ and that the Yaglom limit converges to $\nu$
independently of the initial distribution.  When $\Lambda$ is infinite the
situation is more complex.  Neither existence nor uniqueness of {\qsd} are
guaranteed. An example is the asymmetric random walk $ p=q(i,i+1)= 1-q(i,i-1)$,
for $i\ge 0$. In this case there are infinitely many {\qsd} when $p<1/2$ and
none when $p\ge 1/2$ (see Cavender \cite{cavender} and Ferrari, Martinez and
Picco \cite{fmp} for birth and death more general examples).  For $\Lambda=\N$
under the condition $\lim_{x\to\infty} \P(R>t|Z_0=x) =0$, where $R$ is the
absorption time of $Z_t$, Ferrari, Kesten, Mart\'inez and Picco \cite{fkmp}
prove that the existence of {\qsd} is equivalent to the existence of a positive
exponential moment for $R$, i.e. $\E e^{\theta R}<\infty$ for some $\theta>0$.
When the Yaglom limit exists, it
is known it is a {\qsd}, but existence of the limit is not known in general for
infinite state space.  Phil Pollett maintains an updated bibliography on {\qsd}
in the site http://www.maths.uq.edu.au/\~{}pkp/papers/qsds/qsds.html.

Define the \emph{ergodicity coefficient} of the chain $Q$ by
\begin{equation}
  \label{a7}
  \alpha= \alpha(Q):= \sum_{z\in\Lambda} \inf_{x\in\Lambda\setminus\{z\}} q(x,z)
\end{equation}
If $\alpha(z):=\inf_{x\neq z} q(x,z)>0$, then $z$ is called Doeblin state.  Define the
\emph{maximal absorbing rate} of $Q$ by
\begin{equation}
  \label{a8}
  C=C(Q):=  \sup_{x\in\Lambda} q(x,0)
\end{equation}
Since the chain is absorbed with probability one, $C>0$. On the other hand,
$C<\bar q$, the maximal rate.

Jacka and Roberts \cite{jacka} proved that if there exists a Doeblin state
$z\in\Lambda$ such that $\alpha(z)>C$ and \emph{if} there exists a {\qsd} $\nu$
for $Q$, then $\nu$ is the unique {\qsd} for $Q$ and the Yaglom limit converges
to $\nu$ for any initial measure $\mu$; their proof also works under the weaker
assumption $\alpha > C$. We show that $\alpha > C$ is a \emph{sufficient}
condition for the existence of a {\qsd} for $Q$.
\begin{teo}
  \label{a30}
  If $\alpha>C$ then there exists a unique {\qsd} $\nu$ for $Q$ and the Yaglom
  limit converges to $\nu$ for any initial measure $\mu$.
\end{teo}
The condition $\alpha>C$ is disjoint to the condition $\lim_x\P(R>t|Z_0=x)=0$,
under which \cite{fkmp} show existence of \qsd. On the other hand, $\alpha>0$
implies that $R$ has a positive exponential moment.

\emph{The Fleming-Viot process} ({\fv}).  Let $N$ be a positive integer and
consider a system of $N$ particles evolving on $\Lambda$.  The particles move
independently, each of them governed by the transition rates $Q$ until
absorption. At most one particle is absorbed at any given time.  When a particle
is absorbed, it returns instantaneously to a state in $\Lambda$ chosen with the
empirical distribution of the particles remaining in $\Lambda$.  In other words,
it chooses one of the other particles uniformly and jumps to its position.
Between absorption times the particles move independently governed by $Q$. This
process has been introduced by Fleming and Viot \cite{fv} and studied by Burdzy,
Holyst and March \cite{burdzy}, Grigorescu and Kang \cite{grigorescu} and
L\"obus \cite{lobus} in a Brownian motion setting. The generator acts on
functions $f:\Lambda^{(1,\dots,N)}\to\R$ as follows
\begin{equation}
  \label{a5}
  \cL^N f(\xi) = \sum_{i=1}^N \sum_{y\in\Lambda\setminus\{\xi(i)\}}
  \Bigl[q(\xi(i),y) + q(\xi(i),0)\, \frac{\eta(\xi,y)}{N-1}\Bigr] (f(\xi^{i,y}) - f(\xi))
\end{equation}
where $\xi^{i,y}(j) = y$ for $j=i$
and $\xi^{i,y}(j) = \xi(j)$ otherwise and
\begin{equation}
  \label{a31}
  \eta(\xi,y) := \sum_{i=1}^N \one\{\xi(i)=y\}. 
\end{equation}
We call $\xi_t$ the process in $\Lambda^{(1,\dots,N)}$ with generator \reff{a5}
and $\eta_t = \eta(\xi_t,\cdot)$ the corresponding unlabeled process on
$\N^\Lambda$; $\eta_t(x)$ counts the number of $\xi$ particles in state $x$ at
time $t$. For $\mu$ a measure on $\Lambda$, we denote $\xi^{N,\mu}_t$ the
process starting with independent identically $\mu$-distributed random variables
$(\xi^{N,\mu}_0(i),\,i=1,\dots,N)$; the corresponding variables
$\eta^{N,\mu}_0(x)$ have multinomial law with parameters $N$ and
$(\mu(x),\,x\in\Lambda)$. The profile of the {\fv} process at time $t$ converges
as $N\to\infty$ to the conditioned evolution of the chain $Z_t$:

\begin{teo}\label{teo:lgn}   
  Let $\mu$ be a probability measure on $\Lambda$. Assume that
  $(\xi^{N,\mu}_0(i),\,i=1,\dots,N)$ are iid with law $\mu$. Then, for $t>0$ and
  $x\in\Lambda$,
\begin{eqnarray}
  \lim_{N\to\infty} \E \Bigl(\frac{\eta^{N,\mu}_t(x)}N -
  \varphi^{\mu}_t(x)\Bigr)^2 = 0 \label{lgn1}
\end{eqnarray}  
\end{teo}   
    
The convergence in probability has been proven for Brownian motions in a compact
domain in \cite{burdzy}. Extensions of this result and the process induced in
the boundary have been studied in \cite{grigorescu} and \cite{lobus}.

When $\Lambda$ is finite, the {\fv} process is an irreducible pure-jump Markov
process on a finite state space. Hence it is ergodic (that is, there exists a
unique stationary measure for the process and starting from any measure, the
process converges to the stationary measure). When $\Lambda$ is infinite,
general conditions for ergodicity are not still established.  We prove the
following result
\begin{teo}
  \label{a10}
  If $\alpha > 0$, then for each $N$ the {\fv} process with $N$ particles is
  ergodic.
\end{teo}

Assume $\alpha>0$. Let $\eta^N$ be a random configuration distributed with
the unique invariant measure for the {\fv} process with $N$ particles.
Our next result says that the empiric profile of the invariant measure for the
{\fv} process converges in $L_2$ to the unique {\qsd} for $Q$.
\begin{teo}
\label{t4}  
Assume $\alpha> C$.  Then there exists a probability measure $\nu$ on $\Lambda$
such that for all $x \in \Lambda$,
\begin{eqnarray} 
\label{eq:convinv} 
\lim_{N\to\infty} \E\Bigl( \frac{\eta^{N}(x)}N- \nu(x)\Bigr)^2=0 \label{t42}
\end{eqnarray} 
Furthermore $\nu$ is the unique {\qsd} for $Q$.
\end{teo}   

\paragraph{Sketch of proofs}

The existence part of Theorem \ref{a30} is a corollary of Theorem \ref{t4}. The
rest is a consequence of Jacka and Robert's theorem (stated later as Theorem
\ref{propo:jacka}).

Theorem \ref{a10} is proven by constructing a stationary version of the process
``from the past'' as in perfect simulation. We do it in Section \ref{sec:dqe}. 

Theorems \ref{teo:lgn} and \ref{t4} are both based on the asymptotic
independence of the $\xi$ particles, as $N\to\infty$. Lemma
\ref{u5} later shows that $\varphi_t$ is the unique solution of the Kolmogorov
forward equations
\begin{eqnarray}   
  \frac{d}{dt}\varphi_t^{\mu}(x)= \sum_{y\in\Lambda}
  \varphi_t^{\mu}(y)[q(y,x)+q(y,0)\varphi_t^{\mu}(x)],\qquad x\in\Lambda    \label{a80}
\end{eqnarray}
From a generator computation, taking $f(\xi) = \eta(\xi,x)$ in \reff{a5},
\begin{equation}
  \label{a34}
  \frac{d}{dt} \E\Bigl(\frac{\eta^{N,\mu}_t(x)}N\Bigr)  =  \sum_{y\in\Lambda}
  \E \Bigl(\frac{\eta_t^{N,\mu}(y)}N\Bigl(q(y,x)+q(y,0)\frac{\eta_t^{N,\mu}(x)}{N-1}\Bigr) \Bigr) 
\end{equation}
If solutions of \reff{a34} converge along subsequences as $N\to\infty$, then the
limits equal the unique solution of \reff{a80}. In fact, we prove in Proposition
\ref{corol} that for $x,y\in\Lambda$,
\begin{equation}
  \label{a35}
  \E \Bigl(\eta_t^{N,\mu}(y)\,\eta_t^{N,\mu}(x)-\E\eta_t^{N,\mu}(y)\,\E\eta_t^{N,\mu}(x)\Bigr) = O(N)
\end{equation}
This argument shows the convergence of the means $\E\eta_t^{N,\mu}(x)/N$ to
$\varphi_t^{\mu}(x)$. Since the variances (\reff{a35}, with $x=y$) divided by
$N^2$ go to zero, the $L_2$ convergence follows.

The stationary case is proven analogously.  If $\eta^N$ is distributed with the
invariant measure for the {\fv} process, from \reff{a34},
\begin{equation}
  \label{a38}
  0  =  \sum_{y\in\Lambda}
  \E \Bigl(\frac{\eta^{N}(y)}N\Bigl(q(y,x)+q(y,0)\frac{\eta^{N}(x)}{N-1}\Bigr) \Bigr) 
\end{equation}
Under the hypothesis $\alpha>C$ we show a result for $\eta^N$ analogous to
\reff{a35} to conclude that solutions of \reff{a38} converge to the unique
solution of \reff{eq:quasi2}.


To show that the limits are probability measures it is necessary to show that
the families of measures $(\frac1N \E\eta^{N,\mu}_t,\, N\in\N)$ and $(\frac1N
\E\eta^N,\, N\in\N)$ are tight; we do it in Section~\ref{sec:tight}.

\paragraph{Comments} 
One interesting point of the Fleming-Viot approach is that it permits to show the
existence of a {\qsd} in the $\alpha>C$ case, a new result as far as we know.

Compared with the results for Brownian motion in a bounded region with absorbing
boundary (Burdzy, Holyst and March \cite{burdzy}, Grigorescu and Kang
\cite{grigorescu} and L\"obus \cite{lobus} and other related works), we do not
have trouble with the existence of the {\fv} process, it is immediate here. On
the other hand those works prove the convergence in probability without
computing the correlations. We prove that the fact that the correlations vanish
asymptotically is sufficient to show convergence in probability. For the moment
we are able to show that the correlations vanish for the stationary state under
the hypothesis $\alpha>C$.  

The conditioned distribution $\varphi^\mu_t$ is not necessarily the same as
$\frac1N\E\eta^{N,\mu}_t$, the expected proportion of particles in the {\fv}
process with $N$ particles. This has been proven in Example 2.1 of \cite{burdzy}
for $\Lambda=\{1,2\}$ and $q(1,0)=q(1,2)=q(2,1)=1$. The {\qsd} $\nu$ for this
chain (the unique solution of \reff{eq:quasi2}) is $\nu(1)=(3-\sqrt5)/2$ and
$\nu(2) = (-1+\sqrt5)/2$.  The unlabeled process $\eta^2_t$ with two particles
assumes values in $\{(1,1),\,(2,0),\,(0,2)\}$ with rates
$a((0,2),(1,1))=a((1,1),(0,2))=a((2,0),(1,1))=2$ and $a((1,1),(2,0))=1$.  The
invariant measure for $\eta^2_t$ gives weight $2/5$ to $(1,1)$ and $(0,2)$ and
weight $1/5$ to $(2,0)$. This implies that in equilibrium the mean proportion of
particles in states $1$ and $2$ are $\rho^2(1)=2/5$ and $\rho^2(2)=3/5$
respectively. Our values for $\nu$ and $\rho^2$ do not agree with those of
\cite{burdzy}, but the conclusion is the same: $\nu\neq\rho^2$, which in turns
implies $\frac12\E\eta^{2,\nu}_t\neq \varphi^{\nu}_t=\nu$ for sufficiently large
$t$, as $\frac12\E\eta^{2,\nu}_t$ converges to $\rho^2$ as $t$ grows. More
generally, for rational rates $q$, the equilibrium mean proportions $\rho^N$
have rational components, as they come from the solution of a linear system with
rational coefficients, while those of $\nu$ may be irrational, as $\nu$ is the
solution of a nonlinear system.

To prove tightness we have classified the $\xi$ particles in types. This already
appears in Burdzy, Holyst and March \cite{burdzy} to show the convergence
result. Our application here is somehow simpler. Curiously our tightness proof
needs the same condition ($\alpha>C$) as the vanishing correlations proof.

\section{Construction of {\fv} process}
\label{sec:dqe}
 
In this section we perform the graphic construction of the {\fv} process
$\xi_t^N$.  Recall $C<\infty$ and $\alpha\ge 0$. Recall $\alpha(z)
=\inf_{x\in\Lambda\setminus\{z\}}q(x,z)$. 

For each $i=1,\dots,N$, we define independent stationary marked Poisson processes
(PP's) on $\R$:
\begin{itemize}  
\item Regeneration times. PP rate $\alpha$: $(a^{i}_n)_{n \in\Z}$, with marks
  $(A^{i}_n)_{n \in\Z}$
\item Internal times. PP rate $\bar q-\alpha$: $(b^{i}_n)_{n \in\Z}$, with marks
  $((B^{i}_n(x),\,x\in\Lambda),\,n \in\Z)$
\item Voter times. PP rate $C$: $(c^{i}_n)_{n \in\Z}$, with marks
  $((C^{i}_n,\,(F^i_n(x),\,x\in\Lambda)),\,n \in\Z)$
\end{itemize}

The marks are independent of the PP's and mutually independent. The
denominations will be transparent later. The marginal laws of the marks are:
\begin{itemize}
\item $\P(A^{i}_n = y) = \alpha(y)/\alpha$, $y\in\Lambda$;
\item $\P(B^{i}_n(x)= y)= \displaystyle{\frac{q(x,y)-\alpha(y)}{\bar q-\alpha}}$, $x\in\Lambda$,
  $y\in\Lambda\setminus\{x\}$;\\
 $\P(B^{i}_n(x)= x)=  1- \sum_{y\in\Lambda\setminus\{x\}}\P(B^{i}_n(x)= y)$.
\item $\P(F^{i}_n(x)=1)= \displaystyle{\frac{q(x,0)}{C}}=1-\P(F^{i}_n(x)=0)$,
  $x\in\Lambda$.
\item $\P(C^{i}_n=j)= \displaystyle{\frac{1}{N-1}}$, $j\neq i$.
\end{itemize}

Denote $(\Omega,\mathcal{F}, \P)$ the space where the product of the marked
Poisson processes has been constructed.  Discard the null event corresponding to
two simultaneous events at any given time.

We construct the process in an arbitrary time interval $[s,t]$. Given the mark
configuration $\omega \in \Omega$ we construct
$\xi_{[s,t]}^{N,\xi}(=\xi^{N,\xi}_{[s,t],\omega})$ in the time interval $[s,t]$ as a
function of the Poisson times and its respective marks and the initial
configuration $\xi$ at time $s$.

The relation of this notation with the one in Theorem \ref{teo:lgn} is the
following:
\begin{equation}
  \label{w1}
  \xi_{t}^{N,\mu} = \xi_{[s,t]}^{N,X}
\end{equation}
where $X=(X_1,\dots,X_N)$ is a random vector with iid coordinates, each
distributed according to $\mu$ on $\Lambda$. That is, for any function
$f:\Lambda^N\to\R$,
\begin{equation}
  \label{w2}
  \E f(\xi_{t}^{N,\mu}) = \sum_\xi [\hbox{$\prod_i \mu(\xi(i))$}]\, \E f(\xi_{[s,t]}^{N,\xi}).
\end{equation}

\paragraph{Construction of $\xi_{[s,t]}^{N,\xi}=\xi_{[s,t],\omega}^{N,\xi}$ }\
 
Since for each particle $i$ there are three Poisson processes with rates $C$,
$\alpha$ and $\bar q-\alpha$, the number of events in the interval $[s,t]$ is
Poisson with mean $N(C+\bar q)$. So the events can be ordered from the earliest
to the latest.

At time $s$ the initial configuration is $\xi$. Then, proceed event by event
following the order as follows:

The configuration does not change between Poisson events.

At each regeneration time $a^i_n$ particle $i$ jumps to state $A^i_n$ regardless
the current configuration.

If at the internal time $b^i_n-$ the state of particle $i$ is $x$, then at time
$b^i_n$ particle $i$ jumps to state $B^i_n(x)$ regardless the position of the other
particles.

If at the voter time $c^i_n-$ the state of particle $i$ is $x$ and $F^i_n(x)=1$,
then at time $c^i_n$ particle $i$ jumps to the state of particle $C^i_n$; if
$F^i_n(x)=0$, then particle $i$ does not jump.

The configuration obtained after using all events is $\xi_{[s,t]}^{N,\xi}$.  The
denominations are now clear. At regeneration times a particle jumps to a new
state independently of the current configuration. At voter times a particle
either jumps to the state of another particle chosen at random or does not jump.
At internal times the particle jumps are indifferent to the position of the
other particles.

\begin{lema}
  \label{b14}
  For each $s\in\R$, the process $(\xi_{[s,t]}^{N,\xi}, \,t\ge s)$ is Markov
  with generator \reff{a5} and initial condition $\xi_{[s,s]}^{N,\xi}=\xi$.
\end{lema}

\paragraph{Proof} Follows from the Markov properties of the Poisson processes;
the rate for particle $i$ to jump from $x$ to $y$ is the sum of three terms: (a)
$\alpha \frac{\alpha(y)}{\alpha}$ (the rate of a regeneration event times the
probability that the corresponding mark takes the value $y$), (b) $(\bar q
-\alpha) \frac{q(x,y)-\alpha(y)}{\bar q -\alpha} $ (the maximal rate of internal
events times the probability that the corresponding mark takes the value $y$)
and (c) $C \frac{q(x,0)}{C} \sum_{j\ne i} \one\{\xi(j)=y\}\frac1{N-1} $ (the
maximal absorption rate times the probability the absorption rate from state $x$
divided by the maximal absorption rate times the empirical probability of state
$y$ for the particles different from $i$). The sum of these three rates is the
rate indicated by the generator (the square brackets in \reff{a5} with
$\xi(i)=x$).  \square

\paragraph{Generalized duality}
For each realization of the marked Poisson processes in the interval $[s,t]$ we
construct a set $\Psi^i_\omega[s,t]\subset\{1,\dots,N\}$ corresponding to the
particles involved at time $s$ with the definition of
$\xi_{[s,t],\omega}^{N,\xi}(i)$. We drop the label $\omega$ in the notation.

Initially $\Psi^i[t,t]=\{i\}$ and look backwards in time for the more recent
$i$-Poisson event at some time $\tau$ in the past of $t$ but more recent than
$s$. If $\tau$ is a regeneration event, then we don't need to go further in the
past to know the state of the $i$ particle, so we erase the $i$ particle from
$\Psi^i[\tau-,t]$. If $\tau$ is the voter event $c^i_n$, its $C^i_n$ mark
pointing to particle $j$, say, then we need to know the state of the particle
$i$ at time $\tau-$ to see which $F^{i}_n$ will be used to decide if the $i$
particle effectively takes the value of particle $j$ or not. Hence, we need to
follow backwards particles $i$ and $j$ and we add the $j$ particle to
$\Psi^i[\tau-,t]$. Then continue this procedure starting from each of the
particles in $\Psi^i[\tau-,t]$. The process backwards finishes if $\Psi^i[r,t]$
is empty for some $r$ smaller than $s$ or if we have processed all marks
involving $i$ in the time interval $[s,t]$. More rigorously:

\paragraph{Construction of $\Psi^i [s,t]$} \

We construct $\Psi^i [s,t]$ backwards in time. Changes occur at Poisson events
and $\Psi^i [s,t]$ is constant between two Poisson events. The construction of
$\Psi^i [s,t]$ depends only on the regeneration and voter events. It ignores the
internal events.

Initially $\Psi^i[t,t]=\{i\}$.

Assume $\Psi^i[r',t]$ has been constructed for all $r'\in[\tau',t]$. Let $\tau$
be the time of the latest Poisson event before $\tau'$.

Set $\Psi^i[r'',t]=\Psi^i[\tau',t]$ for all $r''\in(\tau,\tau']$.

If $\tau<s$ stop, we have constructed $\Psi^i[r,t]$ for all $r\in[s,t]$. If not,
proceed as follows.

If $\tau$ is a regeneration event involving particle $j$ (that is, $\tau=a^j_n$
for some $n$), then set $\Psi^i[\tau,t] =\Psi^i[\tau',t] \setminus\{j\}$.

If $\tau$ is a voter event whose mark points to particle $j$ (that is,
$\tau=c^{j'}_n$ for some $j'$ and $n$ and $C^{j'}_n=j$), then set
$\Psi^i[\tau,t] =\Psi^i[\tau',t] \cup\{j\}$.

This ends the iterative step of the construction. 

For a generic Poisson marked event $m$ let $\tim(m)$ be the time it occurs and
$\lab(m)$ its label; for instance $\tim(c^i_n)= c^i_n$, $\lab(c^i_n)=i$. Define
\begin{equation}
  \label{a40o}
  \omega^i[s,t] = \{m\in\omega\,:\, (\lab(m),\tim(m)+) \in
  \{(\Psi^i_\omega[r,t],r),\, r\in[s,t]\},
\end{equation}
the set of marked events in $\omega$ involved in the value of
$\xi_{[s,t],\omega}^{N,\xi}(i)$ and
\begin{equation}
  \label{a40k}
  \xi^i[s,t] = (\xi(j), \,j\in \Psi^i_\omega[s,t]),
\end{equation}
the initial particles involved in the value of $\xi_{[s,t],\omega}^{N,\xi}(i)$.
 
The generalized duality equation is
\begin{eqnarray}
  \xi_{[s,t],\omega}^{N,\xi} (i)= H(\omega^i[s,t], \xi^i[s,t]).
\end{eqnarray}
There is no explicit formula for $H$ but the important point is that for any
real time $s$, $\xi_{[s,t]}^{N,\xi}(i)$ depends only on a \emph{finite} number
of Poisson events contained in $\omega^i[s,t]$ and on the initial state $\xi(j)$
of the particles $j\in \Psi^i_\omega[s,t]$.
The internal marks involved with the
definition of $\xi$ depend on the initial configuration $\xi$ and the evolution
of the process but in any case are bounded by a Poisson random variable with
mean $\bar q|\Psi^i[s,t]|$.

\paragraph{Proof of Theorem \ref{a10}}
If the number of marks in $\omega^i[-\infty,t]$ is finite with probability one,
then the process
\begin{equation}
  \label{a39}
  \xi^N_{t,\omega}(i) 
  = \lim_{s\to-\infty}H(\omega^i[s,t], \xi^i[s,t]),\quad i\in\{1,\dots,N\},\;\;t\in\R
\end{equation}
is well defined with probability one and does not depend on $\xi$. By
construction $(\xi^N_t,\, t\in\R)$ is a stationary Markov process with generator
\reff{a5}.  Since the law at time $t$ does not depend on the initial
configuration $\xi$, the process admits a unique invariant measure, the law of
$\xi^N_t$. See \cite{loss} for more details about this argument.

The number of points in $\omega^i[-\infty,t]$ is finite if and only if for some
finite $s<t$, $\Psi^i[s,t]=\emptyset$. But since there are $3N$ stationary
finite-intensity Poisson processes, with probability one, for almost all
$\omega$ there is an interval $[s(\omega),s(\omega)+1]$ in the past of $t$ such
that there is at least one regeneration mark for all particle $k$ and there are
no voter marks in that interval. We have used here that the regeneration rate
$\alpha>0$. This guarantees that $\Psi^i[s(\omega),t]=\emptyset$. To conclude
notice that if $\Psi^i[s,t]=\emptyset$, then $\Psi^i[s',t]=\emptyset$ for
$s'<s$. \square

\section{Particle correlations in the {\fv} process}\label{sec:cm}
 
In this section we show that the particle-particle correlations in the {\fv}
process with $N$ particles is of the order of $1/N$.
\begin{propo}
\label{corol}
Let $x,y \in \Lambda$.  For all $t >0$
\begin{eqnarray}  
  \Big|\E\( \frac{\eta_t^{N,\mu}(x)\eta_t^{N,\mu}(y)}{N^2}\)
  -\E\(\frac{\eta_t^{N,\mu}(x)}{N}\)\,\E \(\frac{\eta_t^{N,\mu}(y)}{N}\)\Big| 
  \; <\; \frac{1}{N}\, e^{2Ct}
  \label{a40}
\end{eqnarray}  
Assume $ \alpha > C$. Let $\eta^N$ be
  distributed according to the unique invariant measure for the {\fv} process
  with $N$ particles. Then
\begin{eqnarray}   
  \Big|\E \(\frac{\eta^N(x)\eta^N(y)}{N^2}\)
  -\E\(\frac{\eta^N(x)}{N}\)\,\E \(\frac{\eta^N(y)}{N}\)\Big|
  \; <\;  \frac1{N}\,\frac{\alpha}{\alpha-C}
  \label{a41}
\end{eqnarray}   
\end{propo}   

We introduce a 4-fold coupling $(\Psi^ i [s,t],\Psi^ j [s,t],\hat\Psi^ i
[s,t],\hat\Psi^ j [s,t])$ with $\Psi^i[s,t]=\hat\Psi^i[s,t]$ with the property
``$\hat\Psi^j[s,t] \cap \Psi^i[s,t] = \emptyset$ implies $\Psi^j[s,t] =
\hat\Psi^j[s,t]$'' and such that the marginal process
$(\hat\Psi^i[s,t],\hat\Psi^j[s,t])$ have the same law as two independent
processes with the same marginals as $(\Psi^i[s,t],\Psi^j[s,t])$.  The
construction is analogous to the one in Fern\'andez, Ferrari and Garcia
\cite{loss}.

We use two independent families of marked Poisson processes each with the same
law as the Poisson family used in the graphic construction; the marked events
are called {\it red} and {\it green}. We augment the probability space and
continue using $\P$ and $\E$ for the probability and the expectation with
respect to the product space generated by the red and green events. With these
marked events we construct simultaneously the processes $(\Psi^ i [s,t],\Psi^ j
[s,t],\hat\Psi^ i [s,t],\hat\Psi^ j [s,t])$ and a new process $\rI[s,t]$ as
follows.

Initially set ${ \rI [t,t] =0}$, $\hat\Psi^i[t,t]=\Psi^i[t,t]=i$ and
$\hat\Psi^j[t,t]=\Psi^j[t,t]=j$

Go backwards in time as in the construction of $\Psi^i$ in Section \ref{sec:dqe}
proceeding event by event as follows. Assume $\rI [r',t]$, $\hat\Psi^i[r',t]$,
$\Psi^i[r',t]$, $\hat\Psi^j[r',t]$ and $\Psi^j[r',t]$ have been constructed for
all $r'\in[\tau',t]$. Let $\tau$ be the time of the latest Poisson event before
$\tau'$.

If $\rI[\tau',t]=1$ then: (a) if the event is green, use it to update
$\hat\Psi^i[\tau,t]$, $\Psi^i[\tau,t]$ and $\Psi^j[\tau,t]$ only; (b) if the
event is red, use it only to update $\hat\Psi^j[\tau,t]$.

If $\rI[\tau',t] = 0$ then:

(a) if the event is green, then use it to update $\hat\Psi^i[\tau,t]$,
$\Psi^i[\tau,t]$ and $\Psi^j[\tau,t]$. Use it also to update
$\hat\Psi^j[\tau,t]$ only if (after the updating)
$\hat\Psi^j[\tau,t]\cap\hat\Psi^i[\tau,t]=\emptyset $. Otherwise do not update
$\hat\Psi^i[\tau,t]$ and set $\rI[\tau,t]=1$.

(b) if the event is red do not use it to update $\hat\Psi^i[\tau,t]$,
$\Psi^i[\tau,t]$ and $\Psi^j[\tau,t]$. Use it to update $\hat\Psi^j[\tau,t]$
only if after the updating
$\hat\Psi^j[\tau,t]\cap\hat\Psi^i[\tau,t]\neq\emptyset$; in this case set
$\rI[\tau,t] = 1$. Otherwise do not update $\hat\Psi^i[\tau,t]$ and keep
$\rI[\tau,t]=0$.

The processes so constructed satisfy
\begin{enumerate}
\item $\rI[s,t]$ indicates if the hated processes intersect:
\begin{eqnarray} {\rm \rI }[s,t]= \ind \{\hat\Psi^ j [s,t] \cap
    \hat\Psi^i[s,t] \neq \emptyset \}.
  \end{eqnarray}
\item $\Psi^ i [s,t]$ and $\Psi^ j [s,t]$ are constructed using only the {\it
    green} events.
\item $\hat\Psi^i[s,t]$ is also constructed using the  {\it green} events, hence
  it coincides with  $\Psi^i[s,t]$.
\item $\hat\Psi^j[s,t]$ is constructed with a combination of the red and green
  events in such a way that it coincides with $\Psi^j[s,t]$ as long as possible,
  it is independent of $\hat\Psi^i[s,t]$ and has the same marginal distribution
  of $\Psi^j[s,t]$.  
\end{enumerate}

We use the coupling processes to estimate the covariances of
$\xi_{[s,t]}^{N,\mu}$. Call $\omega^j[s,t]$, $\omega^i[s,t]$,
$\hat\omega^j[s,t]$ and $\hat\omega^i[s,t]$ the set of marked events defined
with \reff{a40o} using $\Psi^j[s,t]$, $\Psi^i[s,t]$, $\hat\Psi^j[s,t]$ and
$\hat\Psi^i[s,t]$ respectively. Take two independent random vectors $X$ and $Y$
with the same distribution as in \reff{w1}, that is, iid coordinates with law
$\mu$. Denote the initial particles defined as in \reff{a40k} by $X^j[s,t]$,
$X^i[s,t]$, $\hat X^j[s,t]$ and $\hat Y^i[s,t]$ as function of $\Psi^j[s,t]$,
$\Psi^i[s,t]$, $\hat\Psi^j[s,t]$ and $\hat\Psi^i[s,t]$ respectively. Denote
$\omega^i$ instead of $\omega^i[s,t]$, $X^i$ instead of $X^i[s,t]$, etc.; we
have
\begin{eqnarray}
  \lefteqn{ \P(\xi_{[s,t]}^{N,\mu}(j)=x,\xi_{[s,t]}^{N,\mu}(i)=y)
    -\P(\xi_{[s,t]}^{N,\mu}(j)=x)\P(\xi_{[s,t]}^{N,\mu}(i)=y)} \nonumber\\
  &&=\P(\xi_{[s,t]}^{N,X}(j)=x,\xi_{[s,t]}^{N,X}(i)=y)
  -\P(\xi_{[s,t]}^{N,X}(j)=x)\P(\xi_{[s,t]}^{N,Y}(i)=y) \label{a21}\\
  &&=\E\( \ind \{H(\omega^j,X^j)=x,\,H(\omega^i,X^i)=y)\}- 
  \ind \{H(\hat\omega^j,\hat X^j)=x),\,H(\hat\omega^i,\hat Y^i)=y)\}\)\nonumber
\end{eqnarray}
If $\rI [s,t]=0$ then $\Psi^ j [s',t]=\hat\Psi^ j [s',t]$ and $\Psi^ i
[s',t]=\hat\Psi^ i [s',t]$ for all $s'\in[s,t]$ and the same holds for the
corresponding $\omega$'s. Also, given $\rI [s,t]=0$, $X^j$ and $Y^i$ depend on
disjoint sets of initial particles. This implies that we can couple $X^i$ and
$Y^i$ in such a way that in the event $\rI [s,t]=0$, $X^i=Y^i$. Hence, taking
absolute values in \reff{a21} we get
\begin{eqnarray}\label{a44} 
  |\P(\xi_{[s,t]}^{N,\mu}(j)=x,\xi_{[s,t]}^{N,\mu}(i)=y)
  -\P(\xi_{[s,t]}^{N,\mu}(j)=x)\P(\xi_{[s,t]}^{N,\mu}(i)=y)|  
  \; \leq \; \P({\rm \rI }[s,t]=1).~~  
\end{eqnarray}     

\begin{lema} {\label{lema:bandt}} 
  For $t \geq 0$ and different particles $i,j
  \in\{1,\dots,N\}$
\begin{eqnarray}   
  \P(\rI[s,t] =1 )
  &\le& \frac{1}{N-1}\ \frac{C}{\alpha-C}\ (1 - e^{2(C-\alpha)(t-s)})  \label{a37}
\end{eqnarray}   
\end{lema} 
\prova 
At time $s$ the process $\rI[s,t]$ jumps from $0$ to $1$ at a rate depending on
$\hat\Psi^i [s,t]$ and $\hat\Psi^j[s,t]$ which is bounded above by
\[
\frac{2C}{N-1}\hat\Psi^i [s,t]\hat\Psi^j[s,t]\,\one\{\rI[s,t]=0\}
\]
Dominating the indicator function by one:
\begin{eqnarray}
  \label{a388}
  \P(\rI[s,t]=0\,|\,\mathcal F_{[s,t]}) &\ge& \exp\Big\{-\frac{2C}{N-1}\int_s^t \hat\Psi^i
  [s',t]\hat\Psi^j[s',t]ds'\Big\} 
\end{eqnarray}
where $\mathcal F_{[s,t]}$ is the sigma field generated by $((\hat\Psi^i
[s',t],\,\hat\Psi^j[s',t]),\, s<s'<t)$. From \reff{a388}, using $1-e^{-a}\le a$
and taking expectations,
\begin{equation}
  \label{a42}
  \P(\rI[s,t]=1) \;\le\; \frac{2C}{N-1}\int_s^t \E\hat\Psi^i  [s',t]\,\E\hat\Psi^j[s',t]ds'
\end{equation}
On the other hand, $\hat\Psi^i[s',t]$ is dominated by the position at time $t-s$
of a random walk that grows by one with rate $C$ and decreases by one with rate
$\alpha$. Hence its expectation is bounded above by $e^{(t-s')(C-\alpha)}$.
Substituting this bound in \reff{a42},
\begin{equation}
  \label{a39a}
  \P(\rI[s,t]=1) \le \frac{2C}{N-1}\int_s^te^{2(C-\alpha)(t-s')}ds'
\end{equation}
which gives \reff{a37}. \square

\paragraph{Proof of Proposition \ref{corol}} Defining
\[
\eta^{N,\mu}_{[s,t]}(x) = \sum_{i=1}^N \one\{\xi^{N,\mu}_{[s,t]}=x\}
\]
Then $\eta^{N,\mu}_{[s,t]}$ has the same law as $\eta^{N,\mu}_{t-s}$ and
$\eta^N$ has the same law as $\eta^{N,\mu}_{[-\infty,t]}$. Hence
\begin{eqnarray}\label{var}   
  \E\(\frac{\eta^{N,\mu}_{[s,t]}(x)\,\eta^{N,\mu}_{[s,t]}(y)}{N^2}\)&=& \frac{1}{N^2}   
  \sum_{i=1}^{N}\sum_{j=1}^N
  \,\P(\xi^{N,\mu}_{[s,t]}(i)=x,\xi^{N,\mu}_{[s,t]}(j)=y)  \nonumber \\   
  \frac{\E\eta^{N,\mu}_{[s,t]}(x)\,\E\eta^{N,\mu}_{[s,t]}(y)}{N^2}  &=&  \frac{1}{N^2}   
  \sum_{i=1}^{N}\sum_{j=1}^N
  \,\P(\xi^{N,\mu}_{[s,t]}(i)=x)\P(\xi^{N,\mu}_{[s,t]}(j)=y)  \nonumber 
\end{eqnarray}   
Using this, \reff{a44} and \reff{a37} with $s=0$ and $\alpha=0$ we get \reff{a40}.

If $\alpha>C$, $\eta^{N,\eta}_{[s,t]}$ converges as $s\to-\infty$ to $\eta^N_t$
a configuration distributed with the unique invariant measure, as in Theorem
\ref{a10}, see \reff{a39} for the corresponding statement for $\xi^N_t$. Hence
the left hand side of \reff{a41} is bounded above by $\P(\rI[-\infty,t]=1)$.
Taking $s=-\infty$ in \reff{a37} we get \reff{a41}.  \square

\section{Tightness} \label{sec:tight} \setcounter{equation}{0} 

In this section we prove tightness for the mean densities as probability
measures in  $\Lambda$, indexed by $N$. 
\begin{propo}\label{a51}  
  For all $t>0$, $x \in \Lambda$, $i=1,\dots,N$ and probability 
  $\mu$ on $\Lambda$ it holds
\begin{eqnarray}  
\frac{\E\eta^{N,\mu}_t(x)}{N} &\leq&  e^{Ct}\sum_{z \in \Lambda} \mu(z) P_t(z,x). \label{a50}
\end{eqnarray}
As a consequence the family of measures  $(\E\eta^{N,\mu}_t/N,\,N\in\N)$ is tight.
\end{propo}

Assume $\alpha>0$ and define the probability measure $\mua$ on $\Lambda$ by
\begin{eqnarray}  
  \mua(x) = \frac{\alpha_x}{\alpha},~~ {x \in \Lambda},  \nonumber
\end{eqnarray}  
where $\alpha_x = \inf_z q(z,x)$.  For $z,x \in \Lambda$ define
\begin{eqnarray}  
\label{a66}
R_\lambda(z,x)=\int_0^\infty \lambda e^{-\lambda t} P_t(z,x)dt.  
\end{eqnarray}  
The matrix $R_\lambda$ represents the semigroup $P_t$ evaluated at a random time
$T_\lambda$ exponentially distributed with rate $\lambda$ independent of
$(Z_t)$.  $R_\lambda(z,x)$ is the probability the process $(Z_t^z)$ be in $x$ at
time $T_\lambda$. The matrix $R$ is substochastic: $\sum_{x\in\Lambda}
R_\lambda(z,x)$ is just the probability of non absorption of $(Z_t^z)$ at the
random time $T_\lambda$.
  
\begin{propo}\label{a54}  
  Assume $\alpha>C$ and let $\rho^N(x)$ be the mean proportion of particles in
  state $x$ under the unique invariant measure for the {\fv} process with $N$
  particles. Then for $x \in \Lambda$,
\begin{equation}  
\label{a90}
  \rho^N(x) ~\leq~ \frac{C}{\alpha-C}\, \mua R_{(\alpha-C)}(x)
\end{equation}
As a consequence, the family of measures $( \rho^N,\,N\in\N)$ is tight. 
\end{propo}


\paragraph{Types} To prove the propositions we introduce the concept of types.
We say that particle $i$ is \emph{type} $0$ at time $t$ if it has not been
absorbed in the time interval $[0,t]$. Particles may change type only at
absorption times.  If at absorption time $s$ particle $i$ jumps over particle
$j$ which has type $k$, then at time $s$ particle $i$ changes its type to $k+1$.
Hence, at time $t$ a particle has type $k$ if at its last absorbing time it
jumped over a particle of type $k-1$. We write
\begin{center} 
  $\tip(i,t)$:= type of particle $i$ at time  $t$.
\end{center} 
The marginal law of $\xi_t^{N,\mu}(i)\one\{\tip(i,t)=0\}$ is the law of the
process $Z_t^{\mu}$:
\begin{eqnarray} 
\label{k0} 
\P(\xi_t^{N,\mu}(i)=x,\tip(i,t)=0)=\sum_{z \in
    \Lambda} \mu(z) P_t(z,x).
\end{eqnarray}

\paragraph{Proof of Proposition \ref{a51}} 
Since 
$
\frac{\E\eta^{N,\mu}_t(x)}{N}= \P(\xi_t^{N,\mu}(i)=x),
$
it suffices to show that for $k \geq 0$
\begin{eqnarray}\label{hip}  
  \P(\xi_t^{N,\mu}(i)=x,\tip(i,t)=k) \leq \frac{(Ct)^k}{k!}\sum_{z \in \Lambda} \mu(z) P_t(z,x)  
\end{eqnarray} 
We proceed by induction. By (\ref{k0}) the statement is true for $k=0$. Assume
(\ref{hip}) holds for some $k \geq 0$. We prove it holds for $k+1$.  Time is
partitioned according to the last absorption time $s$ of the $i$th particle.
The absorption occurs at rate bounded above by $C$. The particle jumps at time
$s$ to a particle $j$ with probability $1/(N-1)$, this particle has type $k$ and
state $y$. Then it must go from $y$ to $x$ in the time interval $[s,t]$ without
being absorbed. Using the Markov property, we get:
\begin{eqnarray}  
  \lefteqn{\P(\xi_t^{N,\mu}(i)=x,\tip(i,t)=k+1)}\label{a48}\\  
  &\le&
  \int_0^t \,C\,\frac{1}{N-1}\sum_{j\neq i} \sum_{y \in \Lambda} \,\P(\xi_s^{N,\mu}(j)=
  y,\tip(j,s)=k)\,P_{t-s}(y,x) \,ds.   \label{a100}
\end{eqnarray}  
The symmetry of the particles allows to cancel the sum over $j$ with
$(N-1)^{-1}$. The recursive hypothesis (\ref{hip}) implies that \reff{a100} equals
\begin{eqnarray}  
   &= & \int_0^t C~\frac{(Cs)^k}{k!}\sum_{z  
    \in \Lambda}   \mu(z) \sum_{y \in \Lambda}P_s(z,y) P_{t-s}(y,x) ds \nonumber\\  
  &= &\frac{(Ct)^{k+1}}{(k+1)!}\,\sum_{z \in \Lambda}\mu(z) P_t(z,x). \label{a49}  
\end{eqnarray}  
by Chapman-Kolmogorov. This completes the induction step.  \square

\paragraph{Proof of Proposition \ref{a54}} 
If $\xi^N$ is distributed according to the unique invariant measure for the
{\fv} process then $\rho^N=\P(\xi^N(i)=x)$.  Since $\alpha>0$ we can construct a
version of the stationary process $\xi_s^N$ such that
$\P(\xi^N(i)=x)=\P(\xi_s^N(i)=x)$, $\forall s$.  We analyze the marginal law of
the particle distribution for each type, as in the proof of Proposition
\ref{a51}. Define the types as before, but when a particle meets a regeneration
mark, then the particle type is reset to~0. In the construction, at that time
the state of the particle is chosen with law $\mua$.

Under the hypothesis $\alpha>C$ the process
\[
((\xi^N_t(i), \tip(i,t)),i=1,\dots,N),\,t\in\R)
\] 
is Markovian and can be constructed in a stationary way as $\xi^N_t$. Hence
\begin{eqnarray}
\label{a61}  
A_k(x) := \P(\xi^N_s(i)=x,\tip(i,s)=k)
\end{eqnarray}
does not depend on $s$.

The regeneration
marks follow a Poisson process of rate $\alpha$ and the last regeneration mark
of particle $i$ before time $s$ happened at time $s-T_\alpha^i$, where $T_\alpha^i$
is exponential of rate $\alpha$. Then,
\begin{eqnarray}  
  A_0(x)
  &=&\int_0^\infty \alpha e^{-\alpha t} \sum_{z \in \Lambda} \mua(z) P_t(z,x)dt 
\;=\; \mua R_\alpha (x).  
\label{eq:k0}  
\end{eqnarray}  
A reasoning similar to \reff{a48}-\reff{a100} implies
\begin{eqnarray}
  A_k(x)&\le&
  \int_0^\infty e^{-\alpha t} C \sum_{z\in\Lambda} A_{k-1}(z) P_s(z,x)\,dt.\label{a62}\\ 
  &=&
  \frac{C}\alpha A_{k-1}R_\alpha(x)\;\le\; \(\frac{C}\alpha\)^k \mua
  R_\alpha^{k+1}(x).\label{a63}
 \end{eqnarray}
 We interpret $R^k_\lambda(z,x)$ as the expectation of $P_{\tau_k}(z,x)$, where
 $\tau_k$ is a sum of $k$ independent random variables with exponential
 distribution of rate $\lambda$.  Summing \reff{a62}, and multiplying and
 dividing by $(\alpha-C)$,
\begin{eqnarray}
  \label{a68}
  \P(\xi^N_s(i)= x) \;\le\; \frac{C}{\alpha-C} \sum_{k=0}^\infty
  \(\frac{C}\alpha\)^k\(1-\frac{C}\alpha\) \mua 
   R_\alpha^{k+1}(x)
\end{eqnarray}
The sum can be interpreted as the expectation of $\mua R_\alpha^K$, where $K$ is a
geometric random variable with parameter $p=1-(C/\alpha)$. Since an independent
geometric$(p)$ number of independent exponentials$(\alpha)$ is
exponential$(\alpha p)$, we get
\begin{eqnarray}
  \label{a69}
  \P(\xi^N_s(i)= x) \;\le\; \frac{C}{\alpha-C}\, \mua R_{\alpha-C}(x).\qquad \square 
\end{eqnarray}

\section{Proofs of theorems}\label{sec:teoremas}  
\setcounter{equation}{0} In this section we prove Theorems \ref{a10} and
\ref{t4}. We start deriving the forward equations for $\varphi^\mu_t$ and show
they have a unique solution.
\begin{lema}\label{u5} The Kolmogorov forwards equations for
  $\varphi_t^{\mu}$ are given by 
\begin{eqnarray}   
  \frac{d}{dt}\varphi_t^{\mu}(x)= \sum_{y\in\Lambda}
  \varphi_t^{\mu}(y)[q(y,x)+q(y,0)\varphi_t^{\mu}(x)] \label{a71}
\end{eqnarray}
These equations have a unique solution.   
\end{lema}   
\prova 
The Kolmogorov forward equations for $P_t$ are:
\begin{equation}
\label{eq:kolmogx}   
  \frac{d}{dt}P_t(z,x)= \sum_{y \in\Lambda} P_t(z,y)~ q(y,x) ,\qquad
  z\in\Lambda,\quad x\in  \Lambda\cup\{0\} 
\end{equation} 
Write $\gamma_t=\sum_{z \in \Lambda}\mu(z)P_t(z,0)$ and differentiate \reff{u7}
to get
\begin{eqnarray}   
  \frac{d}{dt}\varphi_t^{\mu}(x) 
  &=& \frac{\sum_{z
      \in \Lambda}\mu(z)\frac{d}{dt}P_t(z,x)}{1-\gamma_t}
  + \frac{(\frac{d}{dt}\gamma_t)}{1-\gamma_t}\cdot\frac{ \sum_{z \in
      \Lambda}\mu(z)P_t(z,x)}{1-\gamma_t}\nonumber \\
  &=&\frac{\sum_{z \in \Lambda}\mu(z)\sum_{y \in\Lambda} 
    P_t(z,y)~q(y,x)}{1-\gamma_t} \nonumber\\ 
  &&\qquad +\;\frac{\sum_{z \in \Lambda} \mu(z)\sum_{y\in\Lambda}
    P_t(z,y)~ q(y,0) }{1-\gamma_t}
  \cdot\frac{\sum_{z \in \Lambda}\mu(z)P_t(z,x)}{1-\gamma_t}
\end{eqnarray}
which equals \reff{a71}.  

To show uniqueness let $\varphi_t$ and $\psi_t$ be two solutions of \reff{a80}
and $\epsilon_t(x) = |\varphi_t(x)-\psi_t(x)|$. Then $\epsilon_t$ satisfies the
inequation
\begin{equation}
  \label{u1}
  \frac{d}{dt} \epsilon_t(y) \le \sum_{z\in\Lambda} \epsilon_t(z) q(z,y) +
  \sum_{z\in\Lambda} |\varphi_t(z)\varphi_t(y)- \psi_t(z)\psi_t(y)|\, q(z,0)
\end{equation}
Bound the modulus with
$\varphi_t(z)\epsilon_t(y)+ \epsilon_t(z)\psi_t(y)$,
sum \reff{u1} in $y$, call $E_t = \sum_{y\in\Lambda} \epsilon_t(y)$ and use
$q(z,y)\le \bar q$ and $q(z,0)\le C$ to get
\begin{equation}
  \label{u3}
  \frac{d}{dt} E_t \le (2\bar q + 2C) E_t   
\end{equation}
This implies $E_t \le E_0\, e^{(2\bar q + 2C)t }$. Since $E_t\ge 0$ and $E_0=0$,
$E_t=0$ for all $t\ge 0$.  \square

\paragraph{Proof of Theorem \ref{teo:lgn}}
We first show convergence of the means
\begin{equation}
  \label{g1}
  \lim_{N\to\infty} \E\Bigl(\eta^{N,\mu}_t(x)\Bigr) = \varphi^{\mu}_t(x).
\end{equation}
  Sum and subtract $\sum_{y
  \in\Lambda}q(y,0) \E(\etatny)\E(\etatnx)$ to \reff{a34} to get
\begin{eqnarray} \label{a82} \E\cL^N\Bigl(\frac{\eta^{N,\mu}_t(x)}{N}\Bigr) &=
  &\sum_{y \in\Lambda} \frac{\E\eta_t^{N,\mu}(y)}{N}\( q(y,x) + q(y,0)
  \frac{\E\eta_t^{N,\mu}(x)}{N}\)\\
  & &\qquad + \sum_{y \in\Lambda}q(y,0)
  \(\E\Bigl(\eta^{N,\mu}_t(y)\eta^{N,\mu}_t(x)\Bigr)
  -\E\eta^{N,\mu}_t(y)\E\eta^{N,\mu}_t(x)\).\label{a81} 
\end{eqnarray}  
By Proposition \ref{a51}, the family $(\E\eta^{N,\mu}_t(x)/N,\,n\in\N)$ is
tight. Call $\rho^\mu_t$ the limit along a convergent subsequence.  Use
$q(x,y)<\bar q$, $q(y,0)<C$ and \reff{a50} to bound the summands in \reff{a82}
and \reff{a81} by $(\bar q+C) e^{Ct} \mu P_t(x)$ and $Ce^{Ct} \mu P_t(x)$,
respectively and use dominated convergence to take the limits inside the sums.
Proposition \ref{corol} implies \reff{a81} converges to zero as $N$ goes to
infinity for any subsequence. We conclude that the limit along a convergent
subsequence must satisfy
\begin{equation}
  \label{a86}
\lim_N  \frac{\E\cL^N\eta^{N,\mu}_t(x)}{N} = \sum_{y\in\Lambda}
  \rho_t^{\mu}(y)[q(y,x)+q(y,0)\rho_t^{\mu}(x)]
\end{equation}
If $F$ is a bounded function
twice continuously differentiable and with uniformly bounded derivatives, then,
\begin{eqnarray}   
  M^{F}_t= F(\eta_t)-F(\eta_0) -\int_0^t \cL F(\eta_s) ds   
\end{eqnarray}   
is a martingale (see Kipnis and Landim (1999), for instance).  We choose
$F(\eta^N_t)=\frac{\eta_t^N(x)}{N}$, $x \in \Lambda$. Since $\E M^F_t=0$,
\begin{eqnarray}   
  \E F(\eta^N_t)=\E F(\eta^N_0)+ \E \int_0^t  \cL^N\frac{\eta^N_s(x)}{N} ds. \label{a89}  
\end{eqnarray}  
From \reff{a86} we conclude that any limit $\rho^\mu_t$ must satisfy
\begin{eqnarray}   
  \rho^\mu_t(x) =\rho^\mu_0(x) + \int_0^t\sum_{y\in\Lambda}
  \rho_s^{\mu}(y)[q(y,x)+q(y,0)\rho_s^{\mu}(x)]dt. \label{a88}  
\end{eqnarray}  
which implies $\rho^\mu_t$ must satisfy \reff{a80}, the forward equations for
$\varphi^\mu_t$. Since there is a unique solution for this equation, the limit
exists and it is $\varphi^\mu_t$.

Taking $y=x$ in \reff{a40}, the variances asymptotically vanish:
\begin{equation}
  \label{y3}
  \lim_{N\to\infty}\frac {\E [\eta^{N,\mu}_t(x)]^2-[\E\eta^{N,\mu}_t(x)]^2
  }{N^2}=0.\qquad
\end{equation}
This concludes the proof. \square

Uniqueness and the Yaglom limit convergence of Theorem \ref{t4} is a consequence
of the next Theorem.

\begin{teo}[Jacka \& Roberts]
\label{propo:jacka} 
If there exists an $x \in \Lambda$ such that $\alpha >C$ and there exists a
{{\qsd}} $\nu$ for $Q$, then $\nu$ is the unique {\qsd} for $Q$ and the Yaglom
limit \reff{yl} converges to $\nu$ for any initial distribution $\mu$.
\end{teo} 

Jacka and Roberts \cite{jacka} use the stronger hypothesis $\inf_{y\in\Lambda}
q(y,x)>C$ for some $x\in\Lambda$ but the proof works under the hypothesis
$\alpha>C$.

\paragraph{Proof of Theorem \ref{t4}}
Since $\alpha>0$, the {\fv} process governed by $Q$ is ergodic by Theorem
\ref{a10}. Call $\eta^N$ a random configuration chosen with the unique invariant
measure. Since $\E\cL^N\eta^N(x) = 0$, summing and subtracting
$\sum_{y\in\Lambda}q(y,0)\rho^N(x)\rho^N(y)$, we get
\begin{eqnarray}
  \lefteqn{0=\sum_{y \in\Lambda}\ron(y)[q(y,x)+q(y,0)\ron(x)]} \nonumber\\   
  & &\qquad+\; \sum_{y
    \in\Lambda}q(y,0)\(\E\Bigl(\etany\etanx\Bigr)-\E\Bigl(\etany\Bigr)\E\Bigl(\etanx\Bigr)\) 
  \label{a91} 
\end{eqnarray}  
which holds for any $N$ and $x \in \Lambda$.  By Proposition \ref{a54},
$(\rho^N,\, N\in\N)$ is tight and by \reff{a90} dominated uniformly in $N$ by a
summable sequence. Bounding $\etanx$ by one, the covariances of \reff{a91} are
bounded by $\rho^N(y)$. Hence we can interchange limit with integral in
\reff{a91} and use \reff{a41} to show that the second term in \reff{a91}
vanishes as $N$ goes to infinity. Then any limit $\rho$ by subsequences must
satisfy the {\qsd} equation \reff{eq:quasi2}. Since by Theorem~\ref{propo:jacka}
the solution is unique, the limit $\lim_N\rho^N$ exists and equals the unique
{\qsd}~$\nu$.  \square

\section{Acknowledgements}
We thank Chris Burdzy for attracting our attention to this problem and for nice
discussions. Se also thank Servet Mart\'\i nez and Pablo Groisman for
discussions. This work is partially supported by FAPESP, CNPq and IM-AGIMB.

 \parbox{1\linewidth}
{
\textsc{Instituto de Matem\'atica e Estat\'istica,\\
Universidade de S\~ao Paulo,\\
Caixa Postal 66281,\\
05311-970 S\~ao Paulo, Brazil}\\
\texttt{pablo@ime.usp.br}, \texttt{nevena@ime.usp.br}\\
\texttt{http://www.ime.usp.br/$\tilde{\,\,\,\,\,}$pablo}, 
\texttt{http://www.ime.usp.br/$\tilde{\,\,\,\,\,}$nevena}
 }

\end{document}